\documentclass[a4paper,12pt]{article}
\usepackage{amssymb}
\usepackage{amsthm}
\usepackage{mathrsfs}
\newtheorem{Lemma}{Lemma}

\newtheorem{Theorem}{Theorem}

\newtheorem{Corollary}{Corollary}

\newtheorem{Conjecture}{Conjecture}
\newtheorem{Question}{Question}

\newcommand{\C}{\mathbb{C}}
\newcommand{\mgnbar}{\overline{{M}}_{g,n}}
\newcommand{\monbar}{\overline{{M}}_{0,n}}
\newcommand{\BB}{\mathscr B}

\title{Positive divisors on quotients of $\monbar$ \\ and the Mori cone of $\mgnbar$}
\author{Claudio Fontanari}
\date{}
\begin{document}
\maketitle

\abstract{We prove that if $m \ge n-3$ then every $S_m$-invariant 
F-nef divisor on the moduli space of stable $n$-pointed curves of 
genus zero is linearly equivalent to an effective combination of 
boundary divisors. As an application, we determine the Mori cone 
of the moduli spaces of stable curves of small genus with few 
marked points.}

\section{Introduction} 

The birational geometry of the moduli space $\mgnbar$ of $n$-pointed 
stable curves of genus $g$ is indeed a fashinating but rather elusive 
subject. In particular, the problem of describing its ample cone has 
attracted very much attention in the last decade (see \cite{Faber:96}, 
\cite{KMK:96}, \cite{GKM:02}, \cite{FarGib:03}, \cite{Gibney:08}, 
\cite{Fontanari:05}, and \cite{Morrison:07} for a comprehensive 
overview).   

Recall that $\mgnbar$ has a natural stratification by topological type, 
the codimension $k$ strata corresponding to curves with at least $k$ 
singular points. 

\begin{Conjecture}\label{fulton} \emph{(\cite{GKM:02} (0.2))} 
A divisor on $\mgnbar$ is ample if and only if it has positive intersection 
with all one-dimensional strata.
\end{Conjecture}

The main result (0.3) of \cite{GKM:02} is that Conjecture \ref{fulton} 
holds if the same claim holds for $\overline{M}_{0,g+n}/S_g$, where $S_g$ 
denotes the symmetric group permuting $g$ marked points.  

\begin{Question}\label{zero} \emph{(\cite{GKM:02} (0.13))}
If a divisor on $\monbar$ has non-negative intersection with all 
one-dimensional strata, does it follow that the divisor is linearly 
equivalent to an effective combination of boundary divisors?
\end{Question}

A positive answer to Question~\ref{zero} would imply Conjecture~\ref{fulton} 
"by an induction which is perhaps the simplest and most telling illustration 
of the power of the inductive structure of the set of all spaces $\mgnbar$" 
(\cite{Morrison:07}). Until now, Question~\ref{zero} has been answered in the affirmative for $n \le 6$, but all available approaches (\cite{GKM:02}, 
\cite{FarGib:03}, and \cite{Fontanari:05}) fail for higher $n$, leading to 
believe that "it seems unlikely (...) even when $n=7$" (\cite{Gibney:08}).  

It is worth stressing that the analogous question for effective divisors has 
a negative answer already for $n=6$ (\cite{Vermeire:02}). On the other 
hand, it is known that every $S_m$-invariant effective divisor on 
$\monbar$ is an effective linear combination of boundary classes 
if $m \ge n-2$ (\cite{Rulla:06}). 

In the same spirit, here we present the following result: 

\begin{Theorem}\label{positive} 
For every integer $n \ge 4$ and $m$ such that $n-3 \le m \le n$, 
every $S_m$-invariant divisor on $\monbar$ intersecting non-negatively 
all one-dimensional strata is linearly equivalent to an effective 
combination of boundary divisors. 
\end{Theorem}

This statement generalizes \cite{FarGib:03}, Proposition~8 (where $n=8$ and 
$m=6$), and our argument simplifies its proof (by showing that case (ii) 
never occurs). As a consequence, we determine the Mori cones of moduli 
spaces of $n$-pointed curves of genus $g$ with small invariants $(n,g)$ 
as follows:

\begin{Corollary}\label{Mori} 
The cone of effective curves of $\mgnbar$ is generated by one-dimensional 
strata for $n=1$, $g \le 9$; $n=2$, $g \le 7$; $n=3$, $g \le 5$.
\end{Corollary}

In the next section we address all essential preliminaries, postponing 
both proofs to the last section. Throughout the paper we work over the 
complex field $\C$.    

This research was partially supported by MiUR (Italy).

\section{The tools}

Let $\delta_S$, $S \subset P := \{1, \ldots, n \}$, and $\psi_i$, 
$i=1, \ldots, n$, be the natural divisor classes on $\monbar$. 

\begin{Lemma}\label{restriction}\emph{(\cite{ArbCor:98}, Lemma~3.3)}
Let $\vartheta: \overline{\mathcal{M}}_{0, A \cup \{ q \}} \longrightarrow 
\overline{\mathcal{M}}_{0, P}$
be the map which associates to any $A \cup \{ q \}$-pointed genus zero 
curve the $P$-pointed genus zero curve obtained by glueing to it a 
fixed $A^c \cup \{ r \}$-pointed genus zero curve via identification 
of $q$ and $r$. Then 
$$
\vartheta^*(\delta_B) = \left\{
\begin{array}{ll}
- \psi_q & \mbox{if $B=A$ or $B = A^c$}\\
\delta_B & \mbox{if $B \subset A$ and $B \ne A$}\\
\delta_{B \setminus A^c \cup \{ q \}} & 
\mbox{if $B \supset A^c$ and $B \ne A^c$}\\
0 & \mbox{otherwise.} 
\end{array}
\right.
$$
\end{Lemma}

According to the standard terminology, a divisor $D$ on $\monbar$ 
is F-nef if it has non-negative intersection with all one-dimensional strata.  

\begin{Lemma}\label{intersection}\emph{(\cite{GKM:02}, Theorem 2.1)}
Let 
$$
D = - \sum_{ \vert S \vert \ge 2} b_S \delta_S + 
\sum_{ \vert S \vert = 1} b_S \psi_S
$$
be a divisor on $\monbar$. Then $D$ is F-nef if and only if 
$$
b_I + b_J + b_K + b_L \ge
b_{I \cup J} + b_{I \cup K} + b_{I \cup L} 
$$
for every partition $I \cup J \cup K \cup L = \{1, 2, \ldots, n \}$. 
\end{Lemma}

If $S_m$ acts on $\monbar$ by permuting the last $m$ marked points, 
then we denote by $(i_R, j_S, k_T, *)$ the partition 
$I \cup J \cup K \cup L = \{1, 2, \ldots, n \}$ with 
$R \cup S \cup T \subseteq \{1, \ldots, n-m \}$, $R \subseteq I$, 
$S \subseteq J$, $T \subseteq K$, and $\vert I \vert = i$, 
$\vert J \vert = j$, $\vert K \vert = k$. We also adopt the 
shorthand notation $i := i_\emptyset$.  

The vector space of $S_m$-invariant divisors on $\monbar$ (up to 
linear equivalence) is generated by the boundary divisors
$$
B^i_T = \sum_{ {\vert S \vert = i} \atop {\{1, \ldots, n-m \} \cap S = T}} \delta_S
$$
with $2 \le i \le n-2$, $T \subseteq \{1, \ldots, n-m \}$, 
and $B^i_T = B^{n-i}_{\{1, \ldots, n-m \} \setminus T}$. 

\begin{Lemma}\label{quotient}\emph{(\cite{Rulla:06}, Corollary 3.2)}
If
\begin{itemize}
\item $\BB_{n,n} = \{B^i_\emptyset: 2 \le i \le \lfloor n/2 \rfloor \}$
\item $\BB_{n,n-1} = \{B^i_T: 2 \le i \le \lfloor n/2 \rfloor \}$
\item $\BB_{n,n-2} = \{B^i_T: 2 \le i \le \lfloor n/2 \rfloor \} 
\setminus \{ B^2_{12} \}$
\item $\BB_{n,n-3} = \{B^i_T: 2 \le i \le \lfloor n/2 \rfloor \} 
\setminus \{ B^2_{12}, B^2_{13}, B^3_{123} \}$
\end{itemize}
then $\BB_{n,m}$ is a basis of the vector space of $S_m$-invariant divisors 
on $\monbar$ for $m \ge n-3$. 
\end{Lemma}

\section{The proofs}

\emph{Proof of Theorem~\ref{positive}.} We are going to show that if $D$ 
is a $S_m$-invariant F-nef divisor on $\monbar$ then 
$$
D = \sum [i]_T B^i_T
$$
where the sum runs over all $B^i_T \in \BB_{n,m}$ as in Lemma~\ref{quotient} 
and every coefficient $[i]_T \ge 0$. We set $[n-i]_{\{1, \ldots,n-m\} \setminus T} 
:= [i]_T$ and $[i] := [i]_\emptyset$.

If $m=n$ and $k \le n-3$ we apply Lemma~\ref{intersection} to the partitions 
$(1,k,i,*)$ for $i=1, \ldots, n-k-2$ and we obtain the inequalities
$$
[k+1]+[k+i]+[i+1]-[k]-[i]-[n-k-i-1] \ge 0.
$$
Summation over $i=1, \ldots, n-k-2$ simplifies to 
$$
(n-k)[k+1] \ge (n-k-2)[k]
$$
and since $[1]=0$ an easy induction implies that $[h] \ge 0$ for every $h \le n-2$. 

If instead $m = n-1$, we consider the partitions $(1,k_{\{1\}},i,*)$ and we get   
$$
[k+1]_{\{1\}}+[k+i]_{\{1\}}+[i+1]-[k]_{\{1\}}-[i]-[n-k-i-1] \ge 0.
$$
Substitution $[n-k-i-1]=[k+i+1]_{\{1\}}$ and summation over $i=1, \ldots, n-k-2$ 
yields  
$$
(n-k)[k+1]_{\{1\}} \ge (n-k-2)[k]_{\{1\}}.
$$
Hence $[h]_{\{1\}} \ge 0$ and $[h] = [n-h]_{\{1\}} \ge 0$ for every $h$. 

Let now $m = n-2$ and $k \le n-3$. The same argument as above applied to the 
partitions $(1,k_{\{1,2\}},i,*)$ for $i=1, \ldots, n-k-2$ implies 
$[h]_{\{1,2\}} \ge 0$ for every $h$ by induction from $[2]_{\{1,2\}}=0$. 
Moreover, $[h] =  [n-h]_{\{1,2\}} \ge 0$ for every $h$. On the other hand,
if we try to address also the coefficients $[h]_{\{\alpha\}}$ with 
$\alpha=1,2$ along the same lines, from the partitions $(1,k_{\{\alpha\}},i,*)$
we obtain
$$
(n-k-1)[k+1]_{\{\alpha\}} + [k+1]_{\{1,2\}} \ge (n-k-2)[k]_{\{\alpha\}}.
$$
In particular, since $[2]_{\{1,2\}}=0$ we deduce $[2]_{\{\alpha\}}=0$. 

\noindent \textbf{Claim.} \emph{For $\alpha=1,2$ we have}
$$
3 [3]_{\{\alpha\}} \ge 3 [2]_{\{\alpha\}} + [4]_{\{\alpha\}}
$$
\emph{and if $k \ge 4$ then}
$$
2 [k]_{\{\alpha\}} \ge \left(\frac{(k-2)(k-3)}{2}-1 \right) [2] + 
(k-2)[2]_{\{\alpha\}} + [k+1]_{\{\alpha\}}. 
$$
The positivity of $[k]_{\{\alpha\}}$ for $k \ge 3$ follows by reverse 
induction from the Claim and the basis step $[n-2]_{\{\alpha\}} = 
[2]_{\{1,2\}\setminus\{\alpha\}} \ge 0$.
In order to check the Claim, we consider the partition 
$((k-1)_{\{\alpha\}},1,1,*)$, which gives 
\begin{equation}\label{one}
2 [k]_{\{\alpha\}} + [2] \ge [k-1]_{\{\alpha\}} + [k+1]_{\{\alpha\}}.
\end{equation}
If $k=3$, it is enough to take into account the partition 
$(1_{\{\alpha\}},1,1,*)$ and subtract from (\ref{one}) the 
corresponding inequality. 
If instead $k \ge 4$, we introduce the sequence of partitions 
$((k-2-i)_{\{\alpha\}},1,1,*)$ for $i=1, \ldots, k-3$, 
providing a weighted sum 
$$
\sum_{i=1}^{k-3} i \left(2 [k-1-i]_{\{\alpha\}} + [2] - [k-2-i]_{\{\alpha\}} - 
[k-i]_{\{\alpha\}}\right) \ge 0
$$ 
which simplifies as
\begin{equation}\label{two}
(k-2) [2]_{\{\alpha\}} + \frac{(k-2)(k-3)}{2}[2] - [k-1]_{\{\alpha\}} \ge 0.
\end{equation} 
By subtracting (\ref{two}) from (\ref{one}) we get the Claim and this completes 
the proof that $[h]_{\{\alpha\}} \ge 0$ for $\alpha =1,2$. 

Finally, we turn to the case $m=n-3$. Notice that the partition 
$(1_{\{1\}}, 1_{\{2\}}, 1_{\{3\}},*)$ yields $[2]_{\{2,3\}} \ge 0$ 
since both $[2]_{\{1,2\}}=0$ and $[2]_{\{1,3\}}=0$ and by repeating 
our standard argument for the partitions $(1,k_{\{1,2,3\}},i,*)$
we obtain $[h]_{\{1,2,3\}}\ge0$ inductively from $[3]_{\{1,2,3\}}=0$. 
In order to check the positivity of all remaining coefficients, 
we reduce ourselves to the previous case $m=n-2$ 
by applying Lemma~\ref{restriction} with $P \setminus A = \{1,2\}$, 
$P \setminus A = \{1,3\}$, and $P \setminus A = \{2,3\}$. Indeed, if $D$ 
is F-nef then both $\vartheta^*(D)$ and $\vartheta^*(D)+[2]_{\{2,3\}} \psi_q$
are F-nef (since $\psi_q$ is ample and $[2]_{\{2,3\}} \ge 0$). Moreover, 
both $\vartheta^*(D)$ on $\overline{M}_{0,P\setminus\{1,2\}\cup\{q\}}$ 
and on $\overline{M}_{0,P\setminus\{1,3\}\cup\{q\}}$ and 
$\vartheta^*(D)+[2]_{\{2,3\}} \psi_q$ on 
$\overline{M}_{0,P\setminus\{2,3\}\cup\{q\}}$
turn out to be expressed in the basis $\BB_{n,n-2}$ as in Lemma~\ref{quotient} 
(indeed, the coefficient of $B^2_{3,q}$ on 
$\overline{M}_{0,P\setminus\{1,2\}\cup\{q\}}$, of $B^2_{2,q}$ on 
$\overline{M}_{0,P\setminus\{1,3\}\cup\{q\}}$, and of $B^2_{1,q}$ on 
$\overline{M}_{0,P\setminus\{2,3\}\cup\{q\}}$ all vanish since 
$[3]_{\{1,2,3\}} = 0$). Hence the previous case applies and the 
proof is over.  

\qed

\emph{Proof of Corollary~\ref{Mori}.} By \cite{GKM:02} (0.3), it is enough 
to prove the same claim for for $\overline{M}_{0,g+n}/S_g$. As pointed out
in \cite{FarGib:03}, Proposition~6, this reduces to check that for all 
boundary restrictions $\nu: \overline{M}_{0,k} \to \overline{M}_{0,g+n}$ 
with $8 \le k \le g+n$ the pull-back of any F-nef divisor is an effective 
combination of boundary classes. In order to do so, for $k=g+n$ we directly 
apply Theorem~\ref{positive}, while if $8 \le k < g+n$ we need also to notice 
that for any boundary restriction $\nu: \overline{M}_{0,g+n-1} \to 
\overline{M}_{0,g+n}$ the pull-back of a F-nef divisor is a F-nef 
divisor on $\overline{M}_{0,g+n-1}/S_{g-2}$.

\qed

\vspace{0.5cm}

\noindent
Claudio Fontanari \newline
Politecnico di Torino \newline
Dipartimento di Matematica \newline
Corso Duca degli Abruzzi 24 \newline
10129 Torino (Italy) \newline
e-mail: claudio.fontanari@polito.it


\begin{thebibliography}{99}

\bibitem{ArbCor:98} E.~Arbarello and M.~Cornalba: Calculating cohomology
groups of mo\-du\-li spaces of curves via algebraic geometry. 
Inst. Hautes \'Etudes Sci. Publ. Math. 88 (1998), 97--127.

\bibitem{Faber:96} C. Faber: Intersection-theoretical computations on 
$\overline {\mathcal{M}}_g$.  Parameter spaces (Warsaw, 1994),  
71--81, Banach Center Publ., 36, Polish Acad. Sci., Warsaw, 1996.

\bibitem{FarGib:03} G.~Farkas and A.~Gibney: 
The Mori cones of moduli spaces of pointed curves of small genus. 
Trans. Amer. Math. Soc. 355 (2003), 1183--1199. 

\bibitem{Fontanari:05} C. Fontanari: A remark on the ample cone of 
$\overline{M}_{g,n}$.  Rend. Sem. Mat. Univ. Politec. Torino 63 
(2005), 9--14. 

\bibitem{Gibney:08} A. Gibney: Numerical criteria for divisors on 
$\overline{M}_g$ to be ample. Comp. Math. (to appear).

\bibitem{GKM:02} A.~Gibney, S.~Keel, and I.~Morrison: 
Towards the ample cone of $\overline{M}_{g,n}$. 
J. Amer. Math. Soc. 15 (2002), 273--294.

\bibitem{KMK:96} S. Keel, J. McKernan: Contractible Extremal Rays on 
$\overline{M}_{0,n}$. \newline arXiv:alg-geom/9607009 (1996).

\bibitem{Morrison:07} I. Morrison: Mori Theory of Moduli Spaces 
of Stable Curves. Projective Press, New York, 2007. \newline
http://www.projectivepress.com/moduli/moristablecurves.pdf

\bibitem{Rulla:06} W. F. Rulla: Effective cones of quotients of 
moduli spaces of stable $n$-pointed curves of genus zero.  
Trans. Amer. Math. Soc.  358  (2006), 3219--3237. 

\bibitem{Vermeire:02} P. Vermeire: A counterexample to Fulton's conjecture on 
$\overline{M}_{0,n}$. J. Algebra  248  (2002), 780--784.

\end{thebibliography}
\end{document}